\documentclass{amsart}
\usepackage{amsfonts}

\theoremstyle{plain}

\newtheorem{definition}{Definition}

\newtheorem{proposition}{Proposition}
\newtheorem{remark}{Remark}

\numberwithin{equation}{section}
\input{tcilatex}

\begin{document}
\title[On Some Discrete Differential Equations ]{On Some Discrete
Differential Equations}
\author{Dejenie A. Lakew}
\address{Virginia State University}
\email{dlakew@vsu.edu}
\urladdr{http://www.vsu.edu}
\date{May 12, 2008}
\subjclass[2000]{ 39A12}

\begin{abstract}
In this short note, we present few results on the use of the discrete
Laplace transform in solving first and second order intial value problems of
discrete differential equations.
\end{abstract}

\maketitle

\ \ \ \ \ \ \ \ \ \ \ \ \ \ \ \ 

In differential equations classes, we teach the Laplace transform as one of
the tools available to find solutions of linear differential equations with
constant coefficients.

\ \ 

What we do in this note is to use a noncontinuous (or discrete) Laplace
transform, that generates sequential solutions which are polynomials in $%
\mathbb{N}
$ or quotients of such polynomials.

I hope readers will find the results very interesting.

\ \ \ 

This work is just a modified version of an unpublished work that I did few
years ago.

\ \ 

Let $f$\ \ be a $%
\mathbb{R}
$-valued sequence $\ f:%
\mathbb{N}
\rightarrow 
\mathbb{R}
$, then

\begin{definition}
The discrete Laplace transform of $f\left( n\right) $ is defined as $\ $%
\begin{equation*}
\ell _{d}\left\{ f\left( n\right) \right\} \left( s\right)
:=\dsum\limits_{n=0}^{\infty }e^{-sn}f\left( n\right) ,\text{ where }s>0.
\end{equation*}

\begin{definition}
The first order difference equation of a sequence $f\left( n\right) $ is
defined as 
\begin{equation*}
\Delta f\left( n\right) :=f\left( n+1\right) -f\left( n\right)
\end{equation*}
\end{definition}
\end{definition}

\begin{proposition}
The first order discrete IVP: 
\begin{equation*}
\triangle f\left( n\right) =n,f\left( 1\right) =1
\end{equation*}
has solution given by

\begin{equation*}
f\left( n\right) =1+\frac{n^{2}-n}{2}.
\end{equation*}
\end{proposition}

\begin{proof}
Taking the transform of both sides of the equation :

$\ell _{d}\left\{ \triangle f\left( n\right) \right\} \left( s\right) =\ell
_{d}\left\{ n\right\} $, we get 
\begin{equation*}
\left( e^{s}-1\right) \ell _{d}\left\{ f\left( n\right) \right\} \left(
s\right) -f\left( 1\right) =\frac{e^{s}}{\left( e^{s}-1\right) ^{2}}.
\end{equation*}

Substituting the value $f\left( 1\right) =1$, and simplifying the expression
we get, 
\begin{equation*}
\ell _{d}\left\{ f\left( n\right) \right\} \left( s\right) =\frac{1}{e^{s}-1}%
+\frac{e^{s}}{\left( e^{s}-1\right) ^{3}}.
\end{equation*}

Then taking the inverse transform we have the solution :

\begin{equation*}
f\left( n\right) =\ell _{d}^{-1}\left\{ \frac{1}{e^{s}-1}+\frac{e^{s}}{%
\left( e^{s}-1\right) ^{3}}\right\}
\end{equation*}

\begin{equation*}
=1+\left( n\ast 1\right) =1+\frac{n^{2}-n}{2}.
\end{equation*}

Here "$\ast $ " is the convolution operator.

\ \ \ \ \ 
\end{proof}

\begin{proposition}
The second order IVP: 
\begin{equation*}
\triangle ^{2}f\left( n\right) =n,f\left( 1\right) =1,\triangle f\left(
1\right) =2,
\end{equation*}

has solution given by :

\begin{equation*}
f\left( n\right) =2n-1+\frac{n\left( n-1\right) \left( n-2\right) }{6}.
\end{equation*}
\end{proposition}

\begin{proof}
First,

\begin{equation*}
\triangle ^{2}f\left( n\right) =\triangle \left( \triangle f\left( n\right)
\right) =f\left( n+2\right) -2f\left( n+1\right) +f\left( n\right)
\end{equation*}

and using the initial conditions we get:

\begin{equation*}
\ell _{d}\left\{ \triangle ^{2}f\left( n\right) \right\} \left( s\right)
=\left( e^{2s}-2e^{s}+1\right) \ell _{d}\left\{ f\left( n\right) \right\}
\left( s\right) -e^{s}-1.
\end{equation*}

\begin{equation*}
\Rightarrow \left( e^{s}-1\right) ^{2}\ell _{d}\left\{ f\left( n\right)
\right\} \left( s\right) -e^{s}-1=\frac{e^{s}}{\left( e^{s}-1\right) ^{2}}
\end{equation*}

\begin{equation*}
\Rightarrow \ell _{d}\left\{ f\left( n\right) \right\} \left( s\right) =%
\frac{1}{\left( e^{s}-1\right) ^{2}}+\frac{e^{s}}{\left( e^{s}-1\right) ^{2}}%
+\frac{e^{s}}{\left( e^{s}-1\right) ^{4}}.
\end{equation*}

\ \ 

Then taking the inverse transform and using convolutions, we get the
solution as :

\begin{equation*}
f\left( n\right) =\left( 1\ast 1\right) +n+\frac{n\left( n-1\right) \left(
n-2\right) }{6}
\end{equation*}

\begin{equation*}
=2n-1+\frac{n\left( n-1\right) \left( n-2\right) }{6}.
\end{equation*}

\ \ \ 
\end{proof}

\begin{proposition}
$\ell _{d}\left\{ \frac{1}{n}\right\} \left( s\right) =s-\ln \left(
e^{s}-1\right) $ \ for $s>0$.
\end{proposition}

\begin{proof}
\begin{equation*}
\ell _{d}\left\{ 1\right\} \left( s\right) =\frac{1}{e^{s}-1}%
=\dsum\limits_{n=1}^{\infty }e^{-sn}\text{ \ for }s>0.
\end{equation*}

Integrating both sides, we get

\begin{equation*}
\ln \left( e^{s}-1\right) -s=\dsum\limits_{n=1}^{\infty }\left( -\frac{e^{sn}%
}{n}\right) =\ell _{d}\left\{ \frac{-1}{n}\right\} \left( s\right)
\end{equation*}

\begin{equation*}
\Rightarrow \ell _{d}\left\{ \frac{1}{n}\right\} \left( s\right) =s-\ln
\left( e^{s}-1\right) .
\end{equation*}

\ \ \ \ 
\end{proof}

We now present discrete IVPs whose solutions are rational sequences in $n$.

\ \ \ \ \ \ \ 

\begin{proposition}
The discrete IVP : 
\begin{equation*}
n\triangle f\left( n\right) =1,f\left( 2\right) =2,\text{for }n\geq 2
\end{equation*}

has solution given by:

\begin{equation*}
f\left( n\right) =1+\dsum\limits_{k=1}^{n-1}\frac{1}{k}.
\end{equation*}
\end{proposition}

\begin{proof}
Taking the transform of both sides of the equation, we have:

\begin{equation*}
\ell _{d}\left\{ n\triangle f\left( n\right) \right\} \left( s\right) =\frac{%
1}{e^{s}-1}.
\end{equation*}

But

\begin{equation*}
\ell _{d}\left\{ n\triangle f\left( n\right) \right\} \left( s\right) =\ell
_{d}\left\{ nf\left( n+1\right) -nf\left( n\right) \right\} \left( s\right)
\end{equation*}

\begin{equation*}
=e^{s}\ell _{d}\left\{ nf\left( n\right) \right\} \left( s\right) -e^{s}\ell
_{d}\left\{ f\left( n\right) \right\} \left( s\right) -\ell _{d}\left\{
nf\left( n\right) \right\} \left( s\right)
\end{equation*}

\begin{equation*}
=\left( e^{s}-1\right) \ell _{d}\left\{ nf\left( n\right) \right\} \left(
s\right) -e^{s}\ell _{d}\left\{ f\left( n\right) \right\} \left( s\right) .
\end{equation*}

Thus,

\begin{equation*}
\left( e^{s}-1\right) \ell _{d}\left\{ nf\left( n\right) \right\} -e^{s}\ell
_{d}\left\{ f\left( n\right) \right\} =\frac{1}{e^{s}-1}.
\end{equation*}

Again,

\begin{equation*}
\ell _{d}\left\{ nf\left( n\right) \right\} \left( s\right) =-\frac{d}{ds}%
\ell _{d}\left\{ f\left( n\right) \right\} \left( s\right) .
\end{equation*}

Therefore,

\begin{equation*}
-\left( e^{s}-1\right) \frac{d}{ds}\ell _{d}\left\{ f\left( n\right)
\right\} \left( s\right) -e^{s}\ell _{d}\left\{ f\left( n\right) \right\}
\left( s\right) =\frac{1}{e^{s}-1},
\end{equation*}

which is an ordinary nonhomogeneous linear differential equation of first
order in $s$ and writing it in standard form we have :

\begin{equation*}
\frac{d}{ds}\ell _{d}\left\{ f\left( n\right) \right\} \left( s\right) +%
\frac{e^{s}}{e^{s}-1}\ell _{d}\left\{ f\left( n\right) \right\} \left(
s\right) =-\frac{1}{\left( e^{s}-1\right) ^{2}}
\end{equation*}

whose solution for $\ell _{d}\left\{ f\left( n\right) \right\} \left(
s\right) $ is given by $\frac{1}{e^{s}-1}-\frac{\ln \left( e^{s}-1\right) -s%
}{e^{s}-1}$.

Then taking the inverse transform , we have :

\begin{equation*}
f\left( n\right) =\ell _{d}^{-1}\left\{ \frac{1}{e^{s}-1}+\frac{s-\ln \left(
e^{s}-1\right) }{e^{s}-1}\right\}
\end{equation*}

\begin{equation*}
=1+\ell _{d}^{-1}\left\{ s-\ln \left( e^{s}-1\right) \right\} \ast \ell
_{d}^{-1}\left\{ \frac{1}{e^{s}-1}\right\}
\end{equation*}

\begin{equation*}
=1+\left( \frac{1}{n}\ast 1\right) =1+\dsum\limits_{k=1}^{n-1}\frac{1}{k}%
\text{ \ for }n>1.
\end{equation*}
\end{proof}

\begin{proposition}
For $n\geq 2$, the IVP : 
\begin{equation*}
\triangle f\left( n\right) =\frac{1}{n^{2}},f\left( 2\right) =2
\end{equation*}

has solution given by

\begin{equation*}
f\left( n\right) =1+\dsum\limits_{k=1}^{n-1}\frac{1}{k^{2}}.
\end{equation*}
\end{proposition}

\begin{proof}
Re-writing the difference equation as : $n\triangle f\left( n\right) =\frac{1%
}{n}$,

taking the transform of both sides and using corollary $3.3$ we get

\begin{equation*}
\frac{d}{ds}\ell _{d}\left\{ f\left( n\right) \right\} \left( s\right) +%
\frac{e^{s}}{e^{s}-1}\ell _{d}\left\{ f\left( n\right) \right\} \left(
s\right) =-\frac{s-\ln \left( e^{s}-1\right) }{e^{s}-1}.
\end{equation*}

Again solving for $\ell _{d}\left\{ f\left( n\right) \right\} \left(
s\right) $, we have

\begin{equation*}
\ \ell _{d}\left\{ f\left( n\right) \right\} \left( s\right) =\frac{1}{%
e^{s}-1}-\frac{1}{e^{s}-1}\int \left( s-\ln \left( e^{s}-1\right) \right) ds.
\end{equation*}

Then solving for $f\left( n\right) $, we have ;

\begin{equation*}
f\left( n\right) =1-\ell _{d}^{-1}\left\{ \frac{1}{e^{s}-1}\right\} \ast
\ell _{d}^{-1}\left\{ \int \left( s-\ln \left( e^{s}-1\right) \right)
ds\right\}
\end{equation*}

\begin{eqnarray*}
&=&1-\left( 1\ast \left( -\frac{1}{n^{2}}\right) \right) \\
&=&1+\dsum\limits_{k=1}^{n-1}\frac{1}{k^{2}}
\end{eqnarray*}

which is the discrete solution of the given discrete differential equation
given in the proposition.
\end{proof}

\begin{remark}
The consequences of the results in this short paper are really far reaching,
if used and expanded further.
\end{remark}

\ \ \ \ \ \ \ \

\end{document}